\RequirePackage[leqno]{amsmath}
\documentclass[leqno]{amsart}
\usepackage{multicol}
\usepackage{enumerate}
\usepackage{color, graphics}
\usepackage{pgfplots}
\usepackage{graphicx, multicol, latexsym, amsmath, amssymb}
\usepackage{grffile}
\usepackage{amsfonts}
 \usepackage[colorlinks=true]{hyperref}
\usepackage{empheq}
\usepackage[T1]{fontenc}
\usepackage{tabularx}
\usepackage[leftcaption]{sidecap}
\sidecaptionvpos{figure}{m}
\usepackage{tikz}
\usepackage{subfigure}
\usepackage[margin=0.5in]{geometry}
\usepackage[pagewise]{lineno}
\usetikzlibrary{backgrounds,automata}
\usepackage{paralist}
\usepackage[inline]{enumitem}

\makeatletter
\@namedef{subjclassname@2020}{\textup{2020} Mathematics Subject Classification}
\makeatother

\def\R{\mathbb R} \def\N{\mathbb N}

\def\R{\mathbb R} \def\N{\mathbb N} 
\def\TM{T_{max}} 
\def
\@cite
#1#2{[{{\bfseries #1}\if@tempswa , #2\fi}]}
\newtheorem{theorem}{Theorem}[section]

\newtheorem{lemma}[theorem]{Lemma}

\newtheorem{definition}[theorem]{Definition}

\title[B\MakeLowercase{oundedness in a chemotaxis model with saturation}] 
      {
\huge{U\MakeLowercase{niform in time} $L^\infty$\MakeLowercase{-estimates for an attraction-repulsion chemotaxis model with double saturation}}}
\author[S. F\MakeLowercase{rassu,} G. V\MakeLowercase{iglialoro and} R. R\MakeLowercase{odr\'iguez} G\MakeLowercase{alv\'an}]{}
\subjclass[2020]{Primary: 35A01, 35Q92. Secondary:  92C17.}
\keywords{Chemotaxis, Attraction-Repulsion, Global existence, Boundedness, Consumption. \\
\textit{$^*$Corresponding author}: giuseppe.viglialoro@unica.it}
\begin{document}
\maketitle

\centerline{S\MakeLowercase{ilvia} F\MakeLowercase{rassu}$^a$, R\MakeLowercase{afael} R\MakeLowercase{odríguez} G\MakeLowercase{alván}$^b$ \and G\MakeLowercase{iuseppe} V\MakeLowercase{iglialoro}$^{a,*}$}
\medskip
{
 \medskip
 \centerline{$^a$Dipartimento di Matematica e Informatica}
 \centerline{Universit\`{a} di Cagliari}
 \centerline{Via Ospedale 72, 09124. Cagliari (Italy)}
 \medskip
 \centerline{$^b$Departamento de Matemáticas}
 \centerline{Universidad de Cádiz}
 \centerline{Campus Universitario de Puerto Real, 11519. 
Cádiz (Spain)}   
}
\bigskip
\begin{abstract}
In this paper we focus on this attraction-repulsion chemotaxis model with consumed signals
\begin{equation}\label{problem_abstract}
\tag{$\Diamond$}
\begin{cases}
u_t= \Delta u - \chi \nabla \cdot (u \nabla v)+\xi \nabla \cdot (u \nabla w)  & \text{ in } \Omega \times (0,T_{max}),\\
v_t=\Delta v- uv  & \text{ in } \Omega \times (0,T_{max}),\\
w_t= \Delta w- uw & \text{ in } \Omega \times (0,T_{max}),
\end{cases}
\end{equation}
formulated in a bounded and smooth domain $\Omega$ of $\R^n$, with $n\geq 2$, for some positive real numbers $\chi,\xi$ and with $\TM\in (0,\infty]$. Once equipped with appropriately smooth initial distributions $u(x,0)=u_0(x)\geq 0$, $v(x,0)=v_0(x)\geq 0$ and $w(x,0)=w_0(x)\geq 0$, as well as Neumann boundary conditions, we establish sufficient assumptions on its data yielding global and bounded classical solutions; these are functions $u,v$ and $w$, with zero normal derivative on $\partial \Omega\times (0,\TM)$, satisfying pointwise the equations in problem \eqref{problem_abstract} with $\TM=\infty$. This is proved for any such initial data, whenever $\chi$ and $\xi$ belong to bounded and open intervals, depending respectively 
	on $\|v_0\|_{L^{\infty}(\Omega)}$ and $\|w_0\|_{L^{\infty}(\Omega)}$. Finally, we illustrate some aspects of the dynamics present within the
chemotaxis system by means of numerical simulations.
%
%
%
%
\end{abstract}
\section{Introduction, main claim and organization of the paper}\label{Intro}
\subsection{The model: discussion and presentation of the main result}\label{IntroDiscussionSection}
Since the advent in the 70's of the biological Keller--Segel models
(\cite{K-S-1970,Keller-1971-MC,Keller-1971-TBC}) idealizing chemotaxis phenomena, many related variants have been finding interest in the mathematical community. Indeed, when facing such problems several technical difficulties appear; of course, this encourages and stimulates researchers in this field.  

In this regard, the model inspiring this article can be formulated as:
\begin{equation}\label{problemAttRep}
\begin{cases}
u_t= \Delta u - \chi \nabla \cdot (u \nabla v)+\xi \nabla \cdot (u
\nabla w)  & \text{ in } \Omega \times (0,T_{max}),\\
\tau v_t= \Delta v -\beta  v +h(u,v)  & \text{ in } \Omega \times (0,T_{max}),\\
\tau w_t= \Delta w-\delta w+k(u,v)& \text{ in } \Omega \times (0,T_{max}),
\end{cases}
\quad \textrm{with} \quad 
\begin{cases}
u_{\nu}=v_{\nu}=w_{\nu}=0 & \text{ on } \partial \Omega \times (0,T_{max}),\\
u(x,0)=u_0(x)& x \in \bar\Omega,\\ v(x,0)=v_0(x) & x \in \bar\Omega, \\  w(x,0)=w_0(x)& x \in \bar\Omega.
\end{cases}
\end{equation}
Herein $\Omega$ is a bounded and smooth domain of $\R^n$, with
$n\in \N$, $\tau \in \{0,1\}$, $\chi, \xi,\beta, \delta>0$ and $h=h(\xi,\eta)$ and
$k=k(\xi,\eta)$ are some regular functions of their arguments $\xi\geq 0$ and $\eta\geq 0$. Moreover, 
further regular initial data  $u_0(x)\geq 0$, $v_0(x)\geq 0$ and $w_0(x)\geq 0$ are as well given, the subscript $\nu$ in $(\cdot)_\nu$ indicates the outward normal derivative on $\partial \Omega$, whereas $T_{max}\in (0,\infty]$ the maximal time up to which solutions to the system are defined. 

In order to properly interpret the mathematical formulation in \eqref{problemAttRep} in the context of real biological phenomena, we specify the following:  if $u=u(x,t)$ denotes  a certain cell distribution (populations, organisms) at the position $x$ and at the time $t$, and $v=v(x,t)$ and $w=w(x,t)$ stand for the concentrations of chemoattractant and chemorepellent (i.e. a chemical signal inducing the cells to, respectively, attract each other and repulse from the other), the system idealizes the motion of the cells, inside an insulated domain (zero-flux on the border: homogeneous Neumann boundary conditions) and initially distributed accordingly to the law of $u_0$. More precisely, such an evolution is influenced by the competition between the aggregation/repulsion impact  (configured at the initial time conforming to  $v_0$ and $w_0$) from the cross terms $\chi u \nabla v/\xi u \nabla w$ both increasing for larger sizes of $\chi$ and $\xi$ and  evolving also following the rates $h(u,v)$ and $k(u,v)$, respectively (second and third equations in \eqref{problemAttRep}). To give an idea of its applicability in the real world,  model \eqref{problemAttRep} has been proposed  to describe the aggregation of microglia involved in the inflammation observed in Alzheimer's disease (\cite{Luca2003Alzheimer}); herein, dimensional, numerical and experimental analyses, in one-dimensional settings, have been proposed for  the fully parabolic version (i.e. $\tau=1$ in both equations) when $h(u,v)$ and $k(u,v)$ are linear functions only of the cell density (formulations of this type are also known as signal-production models).

The aforementioned interplay between the factors describing model \eqref{problemAttRep} strongly influences the cellular motility: in particular, the evolution might relax toward global stabilization and convergence to equilibrium of the cell distribution $u$, but could even degenerate into the so-called \textit{chemotactic collapse}, the mechanism resulting in uncontrolled aggregation processes for $u$, eventually blowing up/exploding at finite time. In particular, moving toward more purely mathematical interpretations (for which we are more interested herein), in the first case, solutions $(u,v,w)$ are defined and bounded for all $(x,t)$ in $\Omega \times (0,\infty)$, in the other, for a certain finite time $\TM$, the solution $(u,v,w)$ ceases to exist for larger value of $\TM$ and becomes unbounded approaching $\TM$. Wanting to provide some details in this framework, for the fully elliptic version (i.e. $\tau=0$ in the equations for $v$ and $w$) of model  \eqref{problemAttRep} 
these results are available in the literature. When linear growths of the
chemoattractant and the chemorepellent are considered, $h(u,v)=\alpha u$, $\alpha>0$,
and $k(u,v)=\gamma u$, $\gamma>0$,
the value $\xi\gamma-\chi\alpha$, measuring 
the difference between the repulsion and attraction contributions, is
critical for $n=2$: particularly, if $\xi\gamma-\chi\alpha>0$
(repulsion dominated regime), in any dimension all solutions
to the model are globally bounded, whereas for
$\xi\gamma-\chi\alpha<0$ (attraction dominated regime) unbounded solutions can be detected (see
\cite{GuoJiangZhengAttr-Rep,LI-LiAttrRepuls,TaoWanM3ASAttrRep,VIGLIALORO-JMAA-BlowUp-Attr-Rep,YUGUOZHENG-Attr-Repul} for some connected studies). On the other hand, for more general production laws, respectively $h$ and $k$ generalizing the prototypes $h(u,v)=\alpha u^s$, $s>0,$ and $k(u,v)=\gamma u^r$, $r\geq 1$, we are only aware of the following recent result, valid for $n\geq 2$ (\cite{ViglialoroMatNacAttr-Repul}): for every  $\alpha,\beta,\gamma,\delta,\chi>0$,  and $r>s\geq 1$ (resp. $s>r\geq 1$), there exists $\xi^*>0$ (resp. $\xi_*>0$) such that if $\xi>\xi^*$ (resp. $\xi\geq \xi_*$), any sufficiently regular initial distribution $u_0(x)\geq 0$ (resp. $u_0(x)\geq 0$ enjoying some smallness assumptions) infers a unique classical and bounded solution. In addition the same conclusion holds true for every  $\alpha,\beta,\gamma,\delta,\chi,\xi>0$, $0<s<1$, $r=1$ and  any sufficiently regular $u_0(x)\geq
0$. 

The chemotactic collapse discussed for the signal-production model \eqref{problemAttRep} is a scenario difficult to detect for alike signal-absorption  models. Indeed, even for the original Keller--Segel system with consumption 
 \begin{equation}\label{problemOriginalKSCosnumption}
\begin{cases}
u_t= \Delta u - \chi \nabla \cdot (u \nabla v)  & \text{ in } \Omega \times (0,T_{max}),\\
v_t=\Delta v-u v & \text{ in } \Omega \times (0,T_{max}),
\end{cases}
\quad \textrm{with} \quad 
\begin{cases}
	u_{\nu}=v_{\nu}=0 & \text{ on } \partial \Omega \times (0,T_{max}),\\
	u(x,0)=u_0(x), \; v(x,0)= v_0(x) & \;x \in \bar\Omega,
\end{cases}
\end{equation}
and corresponding to the simplified two-unknowns version of problem \eqref{problemAttRep} whit $h(u,v)=\beta v-uv$,
the occurrence of blow-up has not been found yet. This is also connected to the observation that comparison arguments suddenly imply  (from the second
equation) that the chemical distribution $v$,  responsible for
gathering phenomena, is uniformly bounded. Despite that, such a bound by itself is not enough to ensure that
classical solutions $(u,v)$ to \eqref{problemOriginalKSCosnumption}
emanating from any sufficiently regular initial data $(u_0,v_0) $ are also
uniformly bounded. Precisely, this holds true only in two-dimensional
settings: this is achieved in \cite{WinklerN-Sto_CPDE}
and \cite{WinklerN-Sto_2d}, where for a more general coupled
chemotaxis-fluid model, respectively globality and convergence of classical solutions are derived. Conversely, in \cite{TaoBoun} it is established that for $n\geq 3$ the
smallness assumption $\chi \lVert v_0\lVert_{L^\infty(\Omega)}\leq
\frac{1}{6(n+1)}$ is required. (This condition has then been improved in  $\chi \lVert v_0\lVert_{L^\infty(\Omega)}<
\frac{\pi}{\sqrt{2(n+1)}}$; see \cite{BaghaeiKhelghatiIranian}.) In particular, attempts to ensure boundedness of solutions to variants tied to problem
\eqref{problemOriginalKSCosnumption} even when $\chi \lVert
v_0\lVert_{L^\infty(\Omega)}$ is larger than $\frac{1}{6(n+1)}$ (or $\frac{\pi}{\sqrt{2(n+1)}}$) have been proposed; we here mention the following two results. In  \cite{LankeitWangConsumptLogistic}, where 
logistic sources with strong damping effects in the equation of the cells have been introduced,  reading
\begin{equation*}
u_t= \Delta u - \chi \nabla \cdot (u \nabla v) +ku-\mu u^2, \quad
\text{ in } \Omega \times (0,T_{max}),\quad k,\mu>0,
\end{equation*}
the authors establish that the resulting Cauchy problem admits classical bounded solutions for
arbitrarily large  $\chi \lVert v_0\lVert_{L^\infty(\Omega)}$ provided
$\mu$ is also larger than a certain expression 
increasing with $\chi \lVert v_0\lVert_{L^\infty(\Omega)}$. But, going towards attraction-repulsion models, when in \eqref{problemAttRep} one considers $h(u,v)\approx \beta v-u^\alpha v$ in the second equation (with $\tau=1$) and  $k(u,v)\approx u^l$ in the third (with $\tau=0$), a model with produced chemorepellent and saturated chemoattractant is obtained;  in \cite{frassuviglialoroConsumptionProduction} boundedness is proved (i) for $l=1$,
$n\in \{1,2\}$, $\alpha\in (0,\frac{1}{2}+\frac{1}{n})\cap (0,1)$ and any $\xi>0$, (ii) for
$l=1$, $n\geq 3$, $\alpha\in (0,\frac{1}{2}+\frac{1}{n})$ and $\xi$ larger than a quantity
depending on $\chi \lVert v_0 \rVert_{L^\infty(\Omega)}$, (iii) for $l>1$ any $\xi>0$, and
in any dimensional settings. 
(Even though herein we are not strictly interested in the influence of logistics in chemotaxis models, we discussed its role in connection to problems of the form in \eqref{problemOriginalKSCosnumption}. Henceforth, we believe that the same is worthwhile to do for formulations quite close to \eqref{problemAttRep}; see, for instance to \cite{ren2019global,ren2020global,ren2020global1}.)

Exactly with the purpose of investigating more deeply the repulsion effect in chemotaxis models like \eqref{problemAttRep}, we therein fix $h(u,v)=\beta v-u v$ and $h(u,v)= \delta w-u w$, thus arriving at this attraction-repulsion system with double saturation:
\begin{equation}\label{problem}
\begin{cases}
u_t= \Delta u - \chi \nabla \cdot (u \nabla v)+\xi \nabla \cdot (u \nabla w)  & \text{ in } \Omega \times (0,T_{max}),\\
v_t= \Delta v - uv  & \text{ in } \Omega \times (0,T_{max}),\\
w_t= \Delta w - uw & \text{ in } \Omega \times (0,T_{max}),\\
u_{\nu}=v_{\nu}=w_{\nu}=0 & \text{ on } \partial \Omega \times (0,T_{max}),\\
u(x,0)=u_0(x), \; v(x,0)=v_0(x), \; w(x,0)=w_0(x) & x \in \bar\Omega.
\end{cases}
\end{equation}
Considering what has been said before, this problem is a natural completion of the above studies on the attraction-repulsion models with only production and  with production and consumption. 

As far as we know,  a general $n$-dimensional theoretical analysis tied to the Cauchy boundary problem \eqref{problem} has not been developed yet;  hereafter, we aim at providing some partial results to extend the corresponding mathematical  comprehension. Specifically, we will derive sufficient conditions involving the data of the model such that the actions of the consumed chemoattractant and chemorepellent guarantee its well-posedness. More precisely, we shall prove this
\begin{theorem}\label{MainTheorem}
For any $n\geq 2$ and some $r>n$, let $\Omega$ be a smooth and bounded domain of $\R^n$ and $(u_0, v_0, w_0)\in (W^{1,r}(\Omega))^3$ any nontrivial initial data  with $u_0, v_0, w_0 \geq 0$ on $\bar{\Omega}$. Then, for $0<\chi <\frac{1}{5n\|v_0\|_{L^{\infty}(\Omega)}}$ and $0<\xi <\frac{1}{5n\|w_0\|_{L^{\infty}(\Omega)}}$ there exists a unique triplet of nonnegative functions
\[
(u, v, w) \in (C^0( [0,\infty); W^{1,r}(\Omega)) \cap C^{2,1}(\bar{\Omega} \times (0,\infty)))^3
\]
solving problem \eqref{problem} and uniformly bounded in time.
\end{theorem}
\subsection{Technical strategy and structure of the article}
In our work we will adapt mathematical reasoning given in \cite{TaoBoun}, and referred to the above introduced consumption chemotaxis model with only attraction.   In both cases, the related steps toward boundedness are connected to some a priori estimates of $\int_\Omega u^k=\int_\Omega u^k(x,t)dx$, for some $k>1$ and some $t>0$, with $u$ being any given solution to problem \eqref{problem} or \eqref{problemOriginalKSCosnumption}. In this direction,  by using boundedness of $v$ and $w$ or $v$ (automatically coming from the parabolic maximum principle),  we can face each problem by suitably employing functionals of the form $\int_\Omega u^k \varphi$. Herein, $\varphi$ is an \textit{ad hoc} function, positive and bounded together with its derivatives; it,  consistently with the investigated model, moreover depends on $(v,w)$ or $v$. Indeed, an evolutive analysis of such a functional leads to a crucial absorption differential inequality in time for the functional itself,  and in turn to the desired bound for $\int_\Omega u^k$. It should be noticed that, evidently, further technical difficulties appear when controlling  
$\frac{d}{dt}\int_\Omega u^k \varphi$ associated to problem \eqref{problem}; in fact it has one more unknown than \eqref{problemOriginalKSCosnumption}.

The rest of the paper is structured as follows. First, in $\S$\ref{LocalSol}  we prove the local existence and uniqueness of a classical solution to model \eqref{problem} and some of its properties. In this same section we give a criterion establishing how to ensure globability and boundedness of local solutions using some of their \textit{a priori} $L^p$-boundedness. In turn, $\S$\ref{EstimatesAndProofSection} is exactly devoted to the derivation of these bounds, by means of which we can deduce the proof of Theorem \ref{MainTheorem}. Finally,  we present some discussions in $\S$\ref{MiscellaneousSection} and, moreover, the theoretical results presented here are investigated numerically in $\S$\ref{SecNumerics}, where two- and three-dimensional simulations are used to discuss how parameters involved in \eqref{problemOriginalKSCosnumption} and \eqref{problem} influence the behaviors of their solutions.

\section{Existence of local-in-time solutions and main properties}\label{LocalSol}
In order to address questions tied to existence and properties of solutions to  problem \eqref{problem}, let us write it in a more suitable form. Specifically,  for $\omega= (\omega_1, \omega_2, \omega_3)=(u,v,w)$ and $\omega(\cdot,0)= (\omega_{1,0}, \omega_{2,0}, \omega_{3,0})=(u_0,v_0,w_0)$, it can be rephrased as 
\begin{equation}\label{Sist}
\begin{cases}
\omega_t = \nabla \cdot (\mathcal{A}(\omega) \nabla \omega) + \mathcal{F}(\omega) &\text{ in } \Omega \times (0, \TM),\\
\omega_{\nu}=0 &\text{ on } \partial \Omega \times (0, \TM),\\
\omega(x, 0)=(u_0,v_0,w_0) & x\in \bar{\Omega},
\end{cases}
\end{equation}
where
\begin{equation*}
\mathcal{A}(\omega)=
\begin{pmatrix}
\mathcal{A}_{11}(\omega) & \mathcal{A}_{12}(\omega) & \mathcal{A}_{13}(\omega)\\
\mathcal{A}_{21}(\omega) & \mathcal{A}_{22}(\omega) & \mathcal{A}_{23}(\omega)\\
\mathcal{A}_{31}(\omega) & \mathcal{A}_{32}(\omega) & \mathcal{A}_{33}(\omega)
\end{pmatrix}=
\begin{pmatrix}
1 & -\chi u & \xi u\\
0 & 1 & 0\\
0 & 0 & 1
\end{pmatrix}
\quad \text{and } \quad
\mathcal{F}(\omega)=
\begin{pmatrix}
\mathcal{F}_1(\omega) \\
\mathcal{F}_2(\omega) \\
\mathcal{F}_3(\omega)
\end{pmatrix} =
\begin{pmatrix}
0 \\
- uv \\
-uw
\end{pmatrix}.
\end{equation*}
With these preparations in our hands, we give the following 
\begin{definition}\label{WeakSolDef}
	For any $n\geq 2$, some $r>n$ and $T>0$, let $\Omega$ be a smooth and bounded domain of $\R^n$ and $(u_0,v_0,w_0) \in (W^{1,r}(\Omega))^3$.
We say that $\omega \in (C([0,T); W^{1,r}(\Omega)))^3$ is a weak $W^{1,r}$-solution of problem \eqref{Sist} in $\Omega \times (0,T)$, if for all test function $\varphi \in C_0^{\infty}(\Omega \times [0,T))$ we have that for $I, J,K = 1, 2, 3$, with $I\neq J \neq K$, 
\[
-\int_0^T \int_\Omega \omega_I \varphi_t - \int_\Omega \omega_{I,0} \varphi(\cdot,0) = 
- \int_0^T \int_\Omega (\mathcal{A}_{II}(\omega) \nabla \omega_I + \mathcal{A}_{IJ}(\omega) \nabla \omega_J + \mathcal{A}_{IK}(\omega) \nabla \omega_K) \cdot \nabla \varphi + \int_0^T \int_\Omega \mathcal{F}_I(\omega) \varphi.
\]
\end{definition}
The forthcoming result ensures that Definition \ref{WeakSolDef} is appropriate in terms of solvability of  problem \eqref{problem}; indeed, it admits weak solutions which are at least local and actually classical. Moreover, they satisfy some uniform bounds. 
\begin{lemma}[\rm{Local existence}]\label{localSol}   
	For any $n\geq 2$ and some $r>n$, let $\Omega$ be a smooth and bounded domain of $\R^n$ and $(u_0, v_0, w_0)\in (W^{1,r}(\Omega))^3$ any nontrivial initial data  with $u_0, v_0, w_0 \geq 0$ on $\bar{\Omega}$.
Then there exists $\TM \in (0,\infty]$ such that problem \eqref{problem} admits a unique nonnegative local-in-time classical solution
\[
(u,v,w)\in (C([0,T_{max}); W^{1,r}(\Omega))\cap  C^{2,1}(\bar{\Omega}\times (0,T_{max})))^3.
\]
In addition, the $u$-component satisfies the mass conservation property, i.e. 
\begin{equation}\label{massConservation}
\int_\Omega u(x, t)dx =\int_\Omega u_0(x)dx=m>0\quad \textrm{for all }\, t \in (0,\TM),
\end{equation}
whilst the components  $v$ and $w$ are such that 
\begin{equation}\label{Cg}
0 \leq v\leq \lVert v_0\rVert_{L^\infty(\Omega)} \quad \textrm{and}\quad 
0 \leq w\leq \lVert w_0\rVert_{L^\infty(\Omega)}\quad \textrm{in}\quad \Omega \times (0,T_{max}).
\end{equation}
Moreover, if there exists some positive constant $C$ such that
\begin{equation}\label{ExtensibilityCriterion}
\|(u(\cdot,t), v(\cdot,t), w(\cdot,t))\|_{L^{\infty}(\Omega)} \leq C \quad \text{for all  } 0<t<\TM,
\end{equation}
then $\TM=\infty$, and $u,v$ and $w$ are uniformly bounded.
\begin{proof}
From the equivalence of the two problems \eqref{problem} and \eqref{Sist}, the existence of a unique maximal weak $W^{1,r}$-solution, in the sense of Definition \ref{WeakSolDef} and
defined in $(0,\TM)$, is consequence of \cite[Theorem 14.4]{AmannBook}. In turn, \cite[Theorem 14.6]{AmannBook} ensures that such a solution is classical. As to its global boundedness, whenever criterion \eqref{ExtensibilityCriterion} is fulfilled, this is deduced by  invoking \cite[Theorem 15.5]{AmannBook}. In addition, since 
$\mathcal{F}_1(0, \omega_2, \omega_3)=0$, $\underline{u} \equiv 0$ is a sub-solution of the first equation in \eqref{problem}, so that from $u_0 \geq 0$ the $u$-component is a super-solution and hence  the parabolic maximum principle implies  $u \geq \underline{u} \equiv 0$ in $\Omega \times (0,\TM)$. Similar arguments provide the bounds in \eqref{Cg}. Finally the mass conservation property follows by integrating over $\Omega$ the first equation of \eqref{problem}, in conjunction with the boundary conditions and the initial data $u_0(x)$.
\end{proof}
\end{lemma}
As it is well known, once local solutions $(u,v,w)$ are given, they are globally bounded if some uniform in time estimate of their components in proper Sobolev spaces 
$W^{k,q}(\Omega)$ can be achieved. Herein, we present for convenience some steps showing that boundedness is implied whenever $u \in L^{\infty}((0,\TM); L^k(\Omega))$
for $k$ larger but arbitrarily close to $\frac{n}{2}$; details and generalizations are available in \cite{BellomoEtAl,HorstWink,W0}.  
\begin{lemma}\label{LocalToUniform}
Under the hypotheses of Lemma \ref{localSol}, let $(u, v, w)$ be such that for some $k>\frac{n}{2}$ the $u$-component belongs to $L^{\infty}((0,\TM); L^k(\Omega))$. Then 
$u,v$ and $w$ are uniformly bounded in time.
\begin{proof}
In view of the criterion \eqref{ExtensibilityCriterion}, it is sufficient to show that  $u,v,w\in L^\infty((0,\TM);L^\infty(\Omega))$ or, by virtue of \eqref{Cg}, solely $u\in L^\infty((0,\TM);L^\infty(\Omega))$. By assuming $\frac{n}{2}<k<n$, we pick $q$ complying with $n<q<k^*:=\frac{nk}{n-k}$. 
By  noting that $-\frac{1}{2}-\frac{n}{2}\left(\frac{1}{k}-\frac{1}{q}\right)> -1$, the representation formula for $v$ provides
\[v(\cdot,t) =e^{t(\Delta-1)}v_0 +\int_0^t e^{(t-s)(\Delta-1)}(1-u(\cdot,s))v(\cdot,s)ds \quad \textrm{for all } \,t\in (0,T_{max}),
\]
and standard properties related to the Neumann heat semigroup $(e^{t \Delta})_{t \geq 0}$ (see Section 2 of \cite{HorstWink} and Lemma 1.3 of \cite{WinklAggre}) lead to 
\begin{equation*}
\begin{split}
\|\nabla v (\cdot, t)\|_{L^q(\Omega)} & \leq \lVert \nabla e^{t (\Delta-1)}v_0 \lVert_{L^q(\Omega)}+\int_0^t \lVert \nabla e^{(t-s)(\Delta-1)}(1-u(\cdot,s))v(\cdot,s)\lVert_{L^q(\Omega)}ds\\ 
&\leq C_S \lVert \nabla v_0\lVert_{L^q(\Omega)}+C_S \lVert v_0 \lVert_{L^{\infty}(\Omega)} \int_0^t (1+(t-s)^{-\frac{1}{2}-\frac{n}{2}(\frac{1}{k}-\frac{1}{q})})e^{-\lambda_1 (t-s)} \lVert 1+u(\cdot,s)\lVert_{L^k(\Omega)}ds \leq c_0. 
\end{split}
\end{equation*}
We have used $v\leq \|v_0\|_{L^{\infty}(\Omega)}$ on $\bar{\Omega}\times (0,\TM)$, taken into account the convergence of  
$\int_0^\infty (t-s)^{-\frac{1}{2}-\frac{n}{2}\big(\frac{1}{k}-\frac{1}{q}\big)}e^{-\lambda_1 (t-s)}ds$ and introduced some $\lambda_1>0$, $C_S>0$ and $c_0>0$. Therefore, we have that $v \in L^{\infty}((0,\TM); W^{1,q}(\Omega))$ and, identically, we also get 
$w \in L^{\infty}((0,\TM); W^{1,q}(\Omega))$. Hence, by setting $\tilde{v}=\chi  v-\xi  w$, the first equation of  \eqref{problem} reads $u_t=\Delta u-\nabla \cdot (u \nabla \tilde{v})$, for any $(x,t)\in \Omega \times (0,T_{max})$, with $ \tilde{v} \in L^\infty((0,T_{max});W^{1,q}(\Omega))$. Through adaptations of \cite[Lemma 3.2]{BellomoEtAl} (see also \cite[Lemma 4.1]{ViglialoroWoolleyAplAnal}), by using some of the above information, we can conclude that $u\in L^\infty((0,\TM);L^\infty(\Omega))$.
\end{proof}
\end{lemma}
\section{A priori estimates and proofs}\label{EstimatesAndProofSection}
With the aim of taking advantage of Lemma \ref{LocalToUniform}, let us dedicate ourselves to the derivation of the specified a priori estimate for the $u$-component. 
\begin{lemma}\label{LemmaEst1}
Let $(u,v,w)$ be the local solution to problem \eqref{problem} provided by Lemma \ref{localSol}, 
$\varphi=\varphi(v,w) \in C^2((0,\|v_0\|_{L^{\infty}(\Omega)})\times(0,\|w_0\|_{L^{\infty}(\Omega)}))$ a function such that $\varphi, \varphi_v, \varphi_w, \varphi_{vv}, \varphi_{ww}, \varphi_{vw}\geq 0$. Then for all positive $\epsilon_1, \epsilon_2, \epsilon_3, \epsilon_4 >0$ the following estimate holds: 
\begin{equation}\label{Est1}
\begin{split}
	&\frac{1}{k}\frac{d}{dt} \int_\Omega u^k \varphi + \left[(k-1)-\epsilon_1^2-\frac{\chi(k-1)\epsilon_2^2}{2}-\epsilon_3^2-\frac{\xi(k-1)\epsilon_4^2}{2}\right] 
	\int_\Omega u^{k-2} \varphi \lvert\nabla u\rvert^2 + \frac{1}{k} \int_\Omega u^k \varphi_{vv} \lvert\nabla v\rvert^2
	\\
	&
	+ \frac{1}{k} \int_\Omega u^k \varphi_{ww} \lvert\nabla w\rvert^2\\
	&\leq \frac{1}{5k} \int_\Omega u^k \varphi_{vv} \lvert\nabla v\rvert^2 \left[\frac{5 k}{\epsilon_1^2} \frac{(\varphi_v)^2}{\varphi \, \varphi_{vv}} 
	+ \frac{5 \chi k(k-1)}{2\epsilon_2^2} \frac{\varphi}{\varphi_{vv}} + 5k \left(\chi + \frac{\xi \|w_0\|_{L^{\infty}(\Omega)}}{8\|v_0\|_{L^{\infty}(\Omega)}}\right) \frac{\varphi_v}{\varphi_{vv}}  
	\right. \\
	&\qquad \qquad \qquad \qquad \quad \left. + \frac{5 \chi k \|w_0\|_{L^{\infty}(\Omega)}}{8\|v_0\|_{L^{\infty}(\Omega)}} \frac{\varphi_w}{\varphi_{vv}}
	+\frac{20 k \|w_0\|_{L^{\infty}(\Omega)}}{\|v_0\|_{L^{\infty}(\Omega)}(k-1)} \frac{\varphi_{vw}}{\varphi_{vv}}\right]\\
	&+ \frac{1}{5k}\int_\Omega u^k \varphi_{ww} \lvert\nabla w\rvert^2 \left[\frac{5 k}{\epsilon_3^2} \frac{(\varphi_w)^2}{\varphi \, \varphi_{ww}} 
	+ \frac{5 \xi k (k-1)}{2\epsilon_4^2)} \frac{\varphi}{\varphi_{ww}} + \frac{10 \chi k \|v_0\|_{L^{\infty}(\Omega)}}{\|w_0\|_{L^{\infty}(\Omega)}} \frac{\varphi_w}{\varphi_{ww}}\right.\\
	&\qquad \qquad \qquad \qquad \qquad \quad \; + \left. \frac{10 \xi k \|v_0\|_{L^{\infty}(\Omega)}}{\|w_0\|_{L^{\infty}(\Omega)}} \frac{\varphi_v}{\varphi_{ww}}
	+\frac{5\|v_0\|_{L^{\infty}(\Omega)}(k-1)}{4k \|w_0\|_{L^{\infty}(\Omega)}}\frac{\varphi_{vw}}{\varphi_{ww}}\right]
	\quad \text{on } (0,\TM).
\end{split}
\end{equation}
\begin{proof}
By a direct calculation and exploiting that $\frac{d}{d t}(u^k\varphi)=ku^{k-1}u_t\varphi + u^k\varphi_v v_t+u^k\varphi_w w_t$, we obtain from problem \eqref{problem} 
\begin{equation}\label{Est2}
\begin{split}
&\frac{1}{k}\frac{d}{dt} \int_\Omega u^k \varphi = \int_\Omega u^{k-1} \varphi \Delta u -\chi \int_\Omega u^{k-1} \varphi \nabla \cdot (u \nabla v)
+\xi \int_\Omega u^{k-1} \varphi \nabla \cdot (u \nabla w) + \frac{1}{k} \int_\Omega u^k \varphi_v \Delta v - \frac{1}{k} \int_\Omega u^{k+1}v \varphi_v\\
&+\frac{1}{k} \int_\Omega u^k \varphi_w \Delta w - \frac{1}{k} \int_\Omega u^{k+1} w \varphi_w\\  
&\leq -(k-1) \int_\Omega u^{k-2} \varphi \lvert \nabla u\rvert^2 - 2 \int_\Omega u^{k-1} \varphi_v \nabla u \nabla v - 2 \int_\Omega u^{k-1} \varphi_w \nabla u \nabla w
+ \chi (k-1) \int_\Omega u^{k-1} \varphi \nabla u \nabla v + \chi \int_\Omega u^k \varphi_v \lvert \nabla v \rvert^2\\ 
&+ \chi \int_\Omega u^k \varphi_w \nabla v \nabla w - \xi (k-1) \int_\Omega u^{k-1} \varphi \nabla u \nabla w - \xi \int_\Omega u^k \varphi_v \nabla v \nabla w 
- \xi \int_\Omega u^k \varphi_w \lvert \nabla w \rvert^2  -\frac{1}{k} \int_\Omega u^k \varphi_{vv} \lvert \nabla v \rvert^2\\ 
&-\frac{2}{k} \int_\Omega u^k \varphi_{vw} \nabla v \nabla w -\frac{1}{k} \int_\Omega u^k \varphi_{ww} \lvert \nabla w \rvert^2 \quad \text{on } (0,\TM). 
\end{split}
\end{equation}
We will now make use of the well-known Young's inequality 
\[C a b \leq \frac{C \epsilon^2}{2} a^2 + \frac{C}{2 \epsilon^2} b^2, \quad \text{valid for } C, a, b >0 \quad \text{and any } \epsilon >0.\]
We have for all $t \in (0,\TM)$
\begin{equation*} 
\begin{split}
- 2 \int_\Omega u^{k-1} \varphi_v \nabla u \nabla v &\leq \epsilon_1^2 \int_\Omega u^{k-2} \varphi \lvert \nabla u\rvert^2 + \frac{1}{\epsilon_1^2} \int_\Omega u^k \frac{(\varphi_v)^2}{\varphi} \lvert \nabla v\rvert^2,\\
\chi (k-1) \int_\Omega u^{k-1} \varphi \nabla u \nabla v &\leq \frac{\chi(k-1)\epsilon_2^2}{2} \int_\Omega u^{k-2} \varphi \lvert \nabla u\rvert^2 + \frac{\chi(k-1)}{2\epsilon_2^2} \int_\Omega u^k \varphi \lvert \nabla v \rvert^2,\\
- 2 \int_\Omega u^{k-1} \varphi_w \nabla u \nabla w &\leq \epsilon_3^2 \int_\Omega u^{k-2} \varphi \lvert \nabla u\rvert^2 + \frac{1}{\epsilon_3^2} \int_\Omega u^k \frac{(\varphi_w)^2}{\varphi} \lvert \nabla w\rvert^2,\\
- \xi (k-1) \int_\Omega u^{k-1} \varphi \nabla u \nabla w &\leq \frac{\xi(k-1)\epsilon_4^2}{2} \int_\Omega u^{k-2} \varphi \lvert \nabla u\rvert^2 + \frac{\xi(k-1)}{2\epsilon_4^2} \int_\Omega u^k \varphi \lvert \nabla w \rvert^2,
\end{split}
\end{equation*}
and
\begin{align*}
\chi \int_\Omega u^k \varphi_w \nabla v \nabla w &\leq \frac{2\chi\|v_0\|_{L^{\infty}(\Omega)}}{\|w_0\|_{L^{\infty}(\Omega)}} \int_\Omega u^k \varphi_{w} \lvert \nabla w\rvert^2 +\frac{\chi \|w_0\|_{L^{\infty}(\Omega)}}{8\|v_0\|_{L^{\infty}(\Omega)}} \int_\Omega u^k \varphi_{w} \lvert \nabla v\rvert^2,\\
- \xi \int_\Omega u^k \varphi_v \nabla v \nabla w &\leq \frac{2\xi\|v_0\|_{L^{\infty}(\Omega)}}{\|w_0\|_{L^{\infty}(\Omega)}} \int_\Omega u^k \varphi_{v} \lvert \nabla w\rvert^2 +\frac{\xi \|w_0\|_{L^{\infty}(\Omega)}}{8 \|v_0\|_{L^{\infty}(\Omega)}} \int_\Omega u^k \varphi_{v} \lvert \nabla v\rvert^2,\\
-\frac{2}{k} \int_\Omega u^k \varphi_{vw} \nabla v \nabla w &\leq \frac{4 \|w_0\|_{L^{\infty}(\Omega)}}{(k-1)\|v_0\|_{L^{\infty}(\Omega)}} \int_\Omega u^k \varphi_{vw} 
\lvert \nabla v\rvert^2 + \frac{(k-1)\|v_0\|_{L^{\infty}(\Omega)}}{4 k^2 \|w_0\|_{L^{\infty}(\Omega)}} \int_\Omega u^k \varphi_{vw} \lvert \nabla w\rvert^2.
\end{align*}
By inserting  all the previous inequalities into \eqref{Est2}, we get for all $t\in (0,\TM)$
\begin{equation*}
\begin{split}
&\frac{1}{k}\frac{d}{dt} \int_\Omega u^k \varphi + \left[(k-1)-\epsilon_1^2-\frac{\chi(k-1)\epsilon_2^2}{2}-\epsilon_3^2-\frac{\xi(k-1)\epsilon_4^2}{2}\right]
\int_\Omega u^{k-2} \varphi \lvert\nabla u\rvert^2 + \frac{1}{k} \int_\Omega u^k \varphi_{vv} \lvert\nabla v\rvert^2+ \frac{1}{k} \int_\Omega u^k \varphi_{ww} \lvert\nabla w\rvert^2 \\
&\leq \int_\Omega u^k \lvert\nabla v\rvert^2 \left[\frac{1}{\epsilon_1^2} \frac{(\varphi_v)^2}{\varphi} + \frac{\chi(k-1)}{2\epsilon_2^2}\varphi + \left(\chi + \frac{\xi \|w_0\|_{L^{\infty}(\Omega)}}{8\|v_0\|_{L^{\infty}(\Omega)}}\right) \varphi_v  + \frac{\chi \|w_0\|_{L^{\infty}(\Omega)}}{8\|v_0\|_{L^{\infty}(\Omega)}} \varphi_w  + \frac{4 \|w_0\|_{L^{\infty}(\Omega)}}{\|v_0\|_{L^{\infty}(\Omega)}(k-1)} \varphi_{vw}\right]\\
&+ \int_\Omega u^k \lvert\nabla w\rvert^2 \left[\frac{1}{\epsilon_3^2} \frac{(\varphi_w)^2}{\varphi} + \frac{\xi(k-1)}{2\epsilon_4^2}\varphi + \frac{2\chi\|v_0\|_{L^{\infty}(\Omega)}}{\|w_0\|_{L^{\infty}(\Omega)}}\varphi_w +\frac{2\xi \|v_0\|_{L^{\infty}(\Omega)}}{\|w_0\|_{L^{\infty}(\Omega)}} \varphi_v + \frac{\|v_0\|_{L^{\infty}(\Omega)}(k-1)}{4k^2\|w_0\|_{L^{\infty}(\Omega)}}\varphi_{vw}\right].
\end{split}
\end{equation*}
\end{proof}
\end{lemma}

\begin{lemma}\label{LemmaEst2}
For some $r>n$, $k>1$ let $u_0, v_0, w_0 \in W^{1,r}(\Omega)$ be nonnegative. If $0<\chi < \frac{1}{10k\|v_0\|_{L^{\infty}(\Omega)}}$ and 
$0<\xi < \frac{1}{10k\|w_0\|_{L^{\infty}(\Omega)}}$, there exists a positive constant $L$ such that 
the local solution $(u,v,w)$ to problem \eqref{problem}, provided by Lemma \ref{localSol}, satisfies
\[
\int_\Omega u^k \leq L \quad \text{on } (0,\TM).
\]
\begin{proof}
Setting 
\begin{equation*}
\varphi(v,w)=e^{\beta^2v^2+\gamma^2w^2}, \quad \text{for } \beta^2 \text{ and } \gamma^2 \text{ to be chosen later},
\end{equation*}
we have 
\begin{equation*}
\varphi_v = 2 \beta^2 v \varphi \quad \varphi_{vv} = 2 \beta^2 \varphi (2 \beta^2 v^2 +1) \quad
\varphi_w = 2 \gamma^2 w \varphi \quad \varphi_{ww} = 2 \gamma^2 \varphi (2 \gamma^2 w^2 + 1)
\quad \varphi_{vw} = \varphi_{wv} = 4 \beta^2 \gamma^2 vw \varphi. 
\end{equation*}
Now, by noting that $\varphi_{vv} \geq 2 \beta^2 \varphi$, we analyze each term of bound \eqref{Est1}, separately,  also in view of \eqref{Cg}. We get: 
\begin{align*}
& I_1 := \frac{5 k}{\epsilon_1^2} \frac{(\varphi_v)^2}{\varphi \, \varphi_{vv}} \leq \frac{10 k}{\epsilon_1^2} \beta^2 \|v_0\|_{L^{\infty}(\Omega)}^2,\\ 
& I_2 := \frac{5 \chi k(k-1)}{2\epsilon_2^2} \frac{\varphi}{\varphi_{vv}} \leq \frac{5\chi k(k-1)}{4\epsilon_2^2\beta^2},\\
& I_3 := 5k \left(\chi + \frac{\xi \|w_0\|_{L^{\infty}(\Omega)}}{8\|v_0\|_{L^{\infty}(\Omega)}}\right) \frac{\varphi_v}{\varphi_{vv}} 
\leq 5k\left(\chi + \frac{\xi \|w_0\|_{L^{\infty}(\Omega)}}{8\|v_0\|_{L^{\infty}(\Omega)}}\right) \|v_0\|_{L^{\infty}(\Omega)},\\
& I_4 := \frac{5 \chi k \|w_0\|_{L^{\infty}(\Omega)}}{8\|v_0\|_{L^{\infty}(\Omega)}} \frac{\varphi_w}{\varphi_{vv}} 
\leq \frac{5 \chi k \|w_0\|_{L^{\infty}(\Omega)}^2}{8\|v_0\|_{L^{\infty}(\Omega)}}\frac{\gamma^2}{\beta^2},\\ 
& I_5 := \frac{20 k \|w_0\|_{L^{\infty}(\Omega)}}{\|v_0\|_{L^{\infty}(\Omega)}(k-1)} \frac{\varphi_{vw}}{\varphi_{vv}} 
\leq \frac{40 k \|w_0\|_{L^{\infty}(\Omega)}^2 \gamma^2}{(k-1)}.
\end{align*}
Similarly, we have $\varphi_{ww} \geq 2 \gamma^2 \varphi$ and  we obtain:
\begin{align*}
& J_1 := \frac{5 k}{\epsilon_3^2} \frac{(\varphi_w)^2}{\varphi \, \varphi_{ww}} \leq \frac{10 k}{\epsilon_3^2} \gamma^2 \|w_0\|_{L^{\infty}(\Omega)}^2,\\ 
& J_2 := \frac{5 \xi k (k-1)}{2\epsilon_4^2} \frac{\varphi}{\varphi_{ww}} \leq \frac{5 \xi k (k-1)}{4\epsilon_4^2 \gamma^2},\\
& J_3 := \frac{10 \chi k \|v_0\|_{L^{\infty}(\Omega)}}{\|w_0\|_{L^{\infty}(\Omega)}} \frac{\varphi_w}{\varphi_{ww}} \leq 10 \chi k \|v_0\|_{L^{\infty}(\Omega)},\\
& J_4 := \frac{10 \xi k \|v_0\|_{L^{\infty}(\Omega)}}{\|w_0\|_{L^{\infty}(\Omega)}} \frac{\varphi_v}{\varphi_{ww}} 
\leq \frac{10 \xi k \|v_0\|_{L^{\infty}(\Omega)}^2}{\|w_0\|_{L^{\infty}(\Omega)}}\frac{\beta^2}{\gamma^2},\\
& J_5 := \frac{5\|v_0\|_{L^{\infty}(\Omega)}(k-1)}{4k \|w_0\|_{L^{\infty}(\Omega)}}\frac{\varphi_{vw}}{\varphi_{ww}} 
\leq \frac{5\|v_0\|_{L^{\infty}(\Omega)}^2(k-1)\beta^2}{2k}.
\end{align*}
Now, let us set 
\begin{equation}\label{BetaGamma}
\beta^2=\frac{\epsilon_1^2}{10k \|v_0\|_{L^{\infty}(\Omega)}^2}, \;
\gamma^2=\frac{\epsilon_3^2}{10k \|w_0\|_{L^{\infty}(\Omega)}^2},
\end{equation}
and 
\begin{equation}\label{ChiXi}
f_1(\epsilon_1,\epsilon_2)=\frac{2 \epsilon_1^2 \epsilon_2^2}{25 k^2(k-1)\|v_0\|_{L^{\infty}(\Omega)}^2} \quad \text{and} \quad
f_2(\epsilon_3,\epsilon_4)=\frac{2 \epsilon_3^2 \epsilon_4^2 }{25 k^2(k-1)\|w_0\|_{L^{\infty}(\Omega)}^2}.
\end{equation}
We herein emphasize that our aim is to figure out which are the largest ranges of $\chi$ and $\xi$ such that, at least in terms of our approach, the validity of this lemma is ensured. Such intervals will result proportional 
to the functions $f_1$ and $f_2$; further,  to an increase of each of these functions corresponds a decrease of the other. Henceforth, in order to simultaneously  optimize both intervals, but at the same time avoiding to make nonpositive the coefficient associated to $\int_\Omega u^{k-2} \varphi \lvert \nabla u \rvert^2$ (crucial in the computations) in estimate \eqref{Est1}, let us consider the open set
\[\mathcal{E}=\left\{(\epsilon_1, \epsilon_2, \epsilon_3, \epsilon_4) \mid \epsilon_1^2+\frac{\chi(k-1)}{2} \epsilon_2^2 + \epsilon_3^2+\frac{\xi(k-1)}{2} \epsilon_4^2 
< k-1\right\}.\]
Evidently, the appropriate choice obeying our requirement entails the supremum of each $\epsilon_1, \epsilon_2, \epsilon_3, \epsilon_4>0$  such that  $\epsilon_1^2=\epsilon_3^2<\frac{k-1}{4}$, $\frac{\chi(k-1)}{2}\epsilon_2^2<\frac{k-1}{4}$ and $\frac{\xi(k-1)}{2}\epsilon_4^2<\frac{k-1}{4}.$ 
In this way, it is seen that if 
\begin{equation}\label{ChiXi1}
0<\chi \leq f_1(\epsilon_1,\epsilon_2) \text{ and }  0<\xi \leq f_2(\epsilon_3,\epsilon_4),  \text{ also reading as } 
0<\chi <\frac{1}{10k\|v_0\|_{L^{\infty}(\Omega)}} \text{ and } 
0<\xi <\frac{1}{10k\|w_0\|_{L^{\infty}(\Omega)}},
\end{equation}
relations \eqref{BetaGamma}, \eqref{ChiXi} and \eqref{ChiXi1} imply that the above terms $I_i$ and $J_i$ comply with
$I_i \leq 1$ and $J_i \leq 1$, with $i=1,\ldots, 5$; in particular  
\begin{equation}
\label{SommaIJ}
\sum_{i=1}^5 I_i \leq 5 \quad \text{and} \quad \sum_{i=1}^5 J_i \leq 5.
\end{equation}
Plugging now the gained estimates \eqref{SommaIJ} into relation \eqref{Est1} yields for some $c_1>0$
\begin{equation}\label{Est5}
\frac{1}{k}\frac{d}{dt} \int_\Omega u^k \varphi + c_1 \int_\Omega u^{k-2} \varphi \lvert\nabla u\rvert^2 \leq 0 \quad \text{for all } t \in (0,\TM). 
\end{equation}
Since $0 \leq v \leq \|v_0\|_{L^{\infty}(\Omega)}$ and $0 \leq w \leq \|w_0\|_{L^{\infty}(\Omega)}$ for all $(x,t)\in \Omega\times (0,\TM)$, we have 
\begin{equation}\label{b}
1\leq \varphi \leq e^{\beta^2 \|v_0\|_{L^{\infty}(\Omega)}^2 + \gamma^2 \|w_0\|_{L^{\infty}(\Omega)}^2}:=b,
\end{equation}
and we obtain 
\begin{equation}\label{Int2}
\int_\Omega u^{k-2} \varphi \lvert\nabla u\rvert^2 \geq \int_\Omega u^{k-2} \lvert\nabla u\rvert^2 = \frac{4}{k^2} \int_\Omega \lvert\nabla u^{\frac{k}{2}}\rvert^2 \quad \text{on }  (0,\TM).
\end{equation}
On the other hand, by applying \eqref{b} and the Gagliardo--Nirenberg inequality (see \cite{Nirenber_GagNir_Ineque}) we deduce that there exists $c_2>0$ such that for
\[
a=\frac{\frac{k}{2}-\frac{1}{2}}{\frac{k}{2}+\frac{1}{n}-\frac{1}{2}} \in (0,1)
\]
we can write
\begin{equation}\label{GN}
\int_\Omega u^k \varphi \leq b \int_\Omega u^k = b \|u^{\frac{k}{2}}\|_{L^2(\Omega)}^2 \leq c_2 \|\nabla u^{\frac{k}{2}}\|_{L^2(\Omega)}^{2a} 
\|u^{\frac{k}{2}}\|_{L^{\frac{2}{k}}(\Omega)}^{2(1-a)} + c_2 \|u^{\frac{k}{2}}\|_{L^{\frac{2}{k}}(\Omega)}^2 \quad \text{on } (0,\TM).
\end{equation}
By exploiting the mass conservation property \eqref{massConservation} and Young's inequality, estimate \eqref{GN} becomes
\begin{equation}\label{GN1}
\int_\Omega u^k \varphi \leq c_3 \left(\int_\Omega \lvert\nabla u^{\frac{k}{2}}\rvert^2\right)^a + c_3 \leq c_4 \int_\Omega \lvert\nabla u^{\frac{k}{2}}\rvert^2 + c_4
\quad \text{for all } t \in (0,\TM)
\end{equation}
and for some positive constants $c_3, c_4$. Therefore, combining bounds \eqref{Int2} and \eqref{GN1}, and introducing proper $c_5, c_6>0$,  lead to
\begin{equation}\label{Int3}
\int_\Omega u^{k-2} \varphi \lvert\nabla u\rvert^2 \geq c_5 \int_\Omega u^k \varphi - c_6 \quad \text{on } (0,\TM).
\end{equation}
By inserting \eqref{Int3} into \eqref{Est5} we have for some $c_7, c_8 >0$
\[
\frac{1}{k}\frac{d}{dt} \int_\Omega u^k \varphi + c_7 \int_\Omega u^k \varphi \leq c_8 \quad \text{for all } t \in (0,\TM),
\]
which, coupled with the initial condition $\frac{1}{k} \int_\Omega u^k(x,0) \varphi = \frac{1}{k} \int_\Omega u_0^k \varphi$, implies $\int_\Omega u^k \varphi \leq c_9$ on $(0,\TM)$ for some $c_9>0$.
Finally, by virtue of expression \eqref{b}, we deduce that for some $L>0$
\begin{equation*}
\int_\Omega u^k \leq \int_\Omega u^k \varphi \leq L \quad \text{on } (0,\TM).
\end{equation*}
\end{proof}
\end{lemma}


\subsubsection*{{\bf{Proof of Theorem \ref{MainTheorem}}}} 
The restrictions  $0<\chi < \frac{1}{5n\|v_0\|_{L^{\infty}(\Omega)}}$ and $0<\xi < \frac{1}{5n\|w_0\|_{L^{\infty}(\Omega)}}$ are such that continuity arguments infer some $k >\frac{n}{2}$ such that the assumption on $\chi$ and $\xi$ of Lemma \ref{LemmaEst2} are complied.  Subsequently, the $u$-component of the local solution $(u,v,w)$ to problem \eqref{problem} belongs to $L^\infty((0,\infty);L^k(\Omega))$, for some $k>\frac{n}{2}$, and Lemma \ref{LocalToUniform} infers the statement. 
\qed
\section{On the use of the functional $\int_\Omega u^k \varphi(v,w)$: miscellaneous}\label{MiscellaneousSection}
The closing part of the paper is concerned with the efficacy of the functional $\int_\Omega u^k \varphi(v,w)$ toward the derivation of the described results.
In this regard we understand that the following comments worth to be analyzed.
\subsection{Improving the condition in \cite[Theorem 1.1]{TaoBoun}} \label{TaoMigliore} 
Since our approach takes ideas from \cite{TaoBoun}, we believe that the following comments may be of interest. The condition $0<\chi \leq \frac{1}{6(n+1)\|v_0\|_{L^{\infty}(\Omega)}}$ in \cite{TaoBoun}, related to the consumption model with only attractive chemical signal, can be improved  to $0<\chi < \frac{2}{3n \|v_0\|_{L^{\infty}(\Omega)}}$. Indeed, assumption \cite[(1.5), Theorem 1.1]{TaoBoun} is derived by taking $k=n+1$ in \cite[Lemma 3.1]{TaoBoun}, but  (see the details in the Lemma \ref{LocalToUniform} above) fixing any $k>\frac{n}{2}$ in \cite[Lemma 3.1]{TaoBoun} it turns in  $0<\chi < \frac{1}{3n \|v_0\|_{L^{\infty}(\Omega)}}$. Even more, this estimate is also more improvable; precisely, for 
\begin{equation}\label{Miso}
\chi \in \left(0,\frac{2}{3n\|v_0\|_{L^{\infty}(\Omega)}}\right)=: \mathcal{I}_{\chi}
\end{equation}
\cite[Theorem 1.1]{TaoBoun} continues valid. 
This can be justified by noting that since in this case the functional $\int_\Omega u^k \varphi(v,w) = \int_\Omega u^k e^{\beta^2v^2+\gamma^2w^2}$ essentially simplifies into $\int_\Omega u^k \varphi(v) = \int_\Omega u^k e^{\beta^2v^2}$, bound \eqref{Est1}  
becomes (see \cite[page 524]{TaoBoun})  on $(0,\TM)$
\begin{equation*}
\begin{split}
&\frac{1}{k}\frac{d}{dt} \int_\Omega u^k \varphi +\left[(k-1)-\tilde{\epsilon}_1^2-\frac{\chi(k-1)\tilde{\epsilon}_2^2}{2}\right]
\int_\Omega u^{k-2} \varphi \lvert\nabla u\rvert^2 
\leq \int_\Omega u^k \lvert\nabla v\rvert^2 \left[\frac{1}{\tilde{\epsilon}_1^2} \frac{(\varphi')^2}{\varphi} + \frac{\chi(k-1)}{2\tilde{\epsilon}_2^2}\varphi + \chi \varphi'-\frac{ \varphi''}{k}\right].
\end{split}
\end{equation*}
Retracing the steps described in Lemma \ref{LemmaEst1}, in order to neglect the term in brackets on the r.h.s. avoiding turning nonpositive the other on the l.h.s., one has to find the supremum of the function  
\[
f_1(\tilde{\epsilon}_1,\tilde{\epsilon}_2) = \frac{2 \tilde{\epsilon}_1^2 \tilde{\epsilon}_2^2}{9k^2(k-1)\|v_0\|_{L^{\infty}(\Omega)}^2} \quad \text{in the open set} \quad 
\tilde{\mathcal{E}}=\left\{(\tilde{\epsilon}_1, \tilde{\epsilon}_2) \bigm| \tilde{\epsilon}_1^2+\frac{\chi(k-1)}{2} \tilde{\epsilon}_2^2 < k-1\right\}.
\]
Such a supremum is exactly $\frac{1}{9k^2\chi\|v_0\|_{L^{\infty}(\Omega)}^2}$, and hence a similar condition as that in \eqref{ChiXi1} entails, for $k>\frac{n}{2}$, the interval defined in \eqref{Miso}.
\subsection{The attraction model vs. the attraction-repulsion model}\label{SubSectionAttrVsAttrRep}
In connection with what we have now done, let us attempt to compare the consumption model with only attraction and its counterpart complemented also with repulsive effects, that is the two models below 
\begin{equation*}
\mathcal{P}_{\chi}:
\begin{cases}
u_t= \Delta u - \chi \nabla \cdot (u \nabla v),\\
v_t= \Delta v - uv, 
\end{cases}
\quad  \text{and}\quad 
\mathcal{P}_{\chi,\xi}:
\begin{cases}
u_t= \Delta u - \chi \nabla \cdot (u \nabla v) + \xi \nabla \cdot (u \nabla w), \\
v_t= \Delta v - uv,\\
w_t= \Delta w - uw, 
\end{cases}
\quad  \text{with}\quad 
\begin{cases}
	\textrm{boundary conditions},\\
	\textrm{initial conditions}.
\end{cases}
\end{equation*}
Using the same notation employed to define the interval in  \eqref{Miso}, let $\mathcal{I}_{\chi,\xi}$ be the interval of $\chi$ where boundedness of solutions to problem 
$\mathcal{P}_{\chi,\xi}$ may be ensured; $\mathcal{I}_{\chi}$, naturally, stands for that related to problem 
$\mathcal{P}_{\chi}$. Even though by biological intuitions it should be conceivable that  $\mathcal{I}_{\chi,\xi}\supset\mathcal{I}_{\chi}$ for every $\xi>0$, this cannot be directly deduced by employing the functionals $\int_\Omega u^k \varphi$. Indeed, by observing from \eqref{ChiXi} that taking 
$(\epsilon_3, \epsilon_4) \to (0,0)$ and $(\epsilon_1,\epsilon_2) \to \left(\sqrt{\frac{k-1}{2}},\sqrt{\frac{1}{\chi}}\right)$ corresponds to reduce/increase as much as possible the repulsion/attraction strength, the conclusion of Theorem \ref{MainTheorem} would apply for positive $\xi$ arbitrarily small, and $\chi$ positive up (depending on $\xi$) to the supremum  value $\frac{2}{5n \lVert v_0\rVert_{L^{\infty}(\Omega)}}$. Nevertheless, the resulting limit interval 
$\left(0,\frac{2}{5n \lVert v_0\rVert_{L^{\infty}(\Omega)}}\right)$ would remain strictly included in $\mathcal{I}_{\chi}$. We understand that the reason of this undesired outcome is connected to the quantity of terms inherited from the functional in terms of the unknowns and equations of the problem: two in $\mathcal{P}_{\chi}$, three in $\mathcal{P}_{\chi,\xi}$. As analyzed, the exceeding part from the larger system produces extra addends implying that certain conclusions hold true for ranges of $\chi$ smaller than those linked to the other system. However, aided by the next simulations, and in particular those explained in $\S$\ref{SubSectionModelDoubleConsSimNumer}, we can conclude that the aforementioned intuition is rather founded.   
\section{Numerical simulations}\label{SecNumerics}
 In this section we numerically test the presented theoretical results by simulating systems \eqref{problemOriginalKSCosnumption} and \eqref{problem} in two and three dimensions (also abbreviated with $2d$ and $3d$). Our computational and quantitative experiments are meant to support and complement the qualitative analyses in the earlier sections.
\subsection{Methodology, difficulties and computational devices}\label{sec:numerical-tests}
Numerical simulations of chemotaxis-type systems, like those now mentioned, are by no means an easy
goal but, on the contrary, a challenging task. Indeed, the solutions to such systems may exhibit interesting mathematical properties which cannot be easily reproduced in a discrete framework. In particular, it is rather difficult to develop numerical schemes preserving desired properties of solutions, like non-negativity or energy laws; moreover, it is also far from being a standard task capturing possible blow-up phenomena. (Some hints describing the technical details behind this complexity may be found, for instance, in  \cite{chertock_second-order_2008,zhou_finite_2017,saad_mazen_finite_2014,gutierrez-santacreu_analysis_2021,epshteyn_new_2009}.) 

In our forthcoming $2d$/$3d$ simulations for both models \eqref{problem} and \eqref{problemOriginalKSCosnumption}, we introduce a Finite Element Method (FEM) space discretization, defined on a triangular/tetrahedral mesh with size $h$ of the domain $\Omega$, with first order piecewise polynomials. On the other hand, we  consider an implicit Euler scheme for the time discretization $\Delta t$, where each time step is split into two stages. First we compute $v$ and $w$ as solutions of the (uncoupled) linear parabolic equations modeling the consumption, being $u$ explicitly known from the previous time step. Then the values of $u$ at any node of the mesh are updated by solving the chemotaxis equation, exactly using the gained quantities of the attraction and repulsion terms $v$ and $w$ at those nodes. (We specify that for the computation we used the open source library FEniCS (\cite{LangtangenLogg2017Fenics}) and a parallel solver has been employed at each time step for the resolution of the FEM linear systems.) It is  precisely in this latter chemotaxis stage where the more severe difficulties arise. In fact, the numerical treatment of the taxis-driven terms may conduce to undesirable effects (especially in three-dimensional settings), as the mentioned loss of positivity of $u$ and/or related spurious oscillations, with consequent disastrous results. However, the combination of the Flux-Corrected Transport technique (\cite{strehl_flux-corrected_2010}) and the use of proper devices (in the specific, a computer with two processors AMD Opteron 6320, 32 cores and  1Tb of RAM) allow us to circumvent these situations, hence achieving qualitatively acceptable results. 
\subsection{Simulations of model \eqref{problemOriginalKSCosnumption}}
Through the following simulations, we essentially aim at corroborating the absence of blow-up solutions to the classical consumption model even for large values of the chemoattractant $\chi$ associated to the attractive chemical agent $v$.
\subsubsection{Two-dimensional numerical tests.}
\label{sec:2d-numerical-test}
\begin{figure}
\centering
\begin{tabular}[c]{cc}
\rotatebox{90}{\quad $t=0$}& \includegraphics[width=0.73\linewidth]{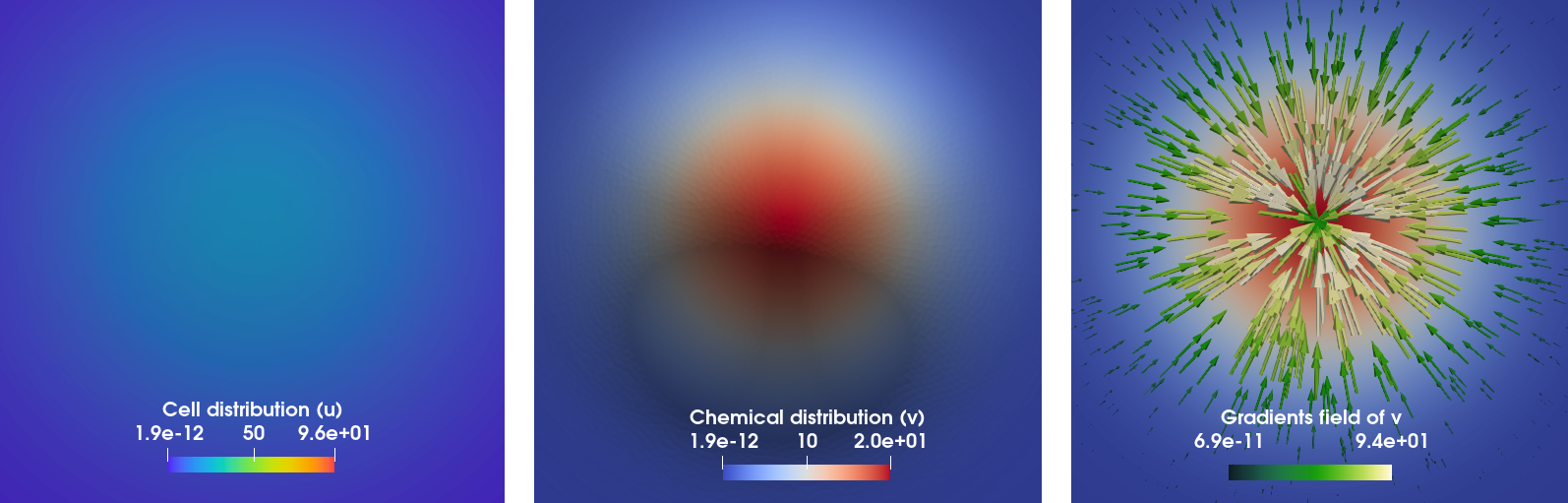}
  \\
\rotatebox{90}{\quad $t=0.0001$}&\includegraphics[width=0.73\linewidth]{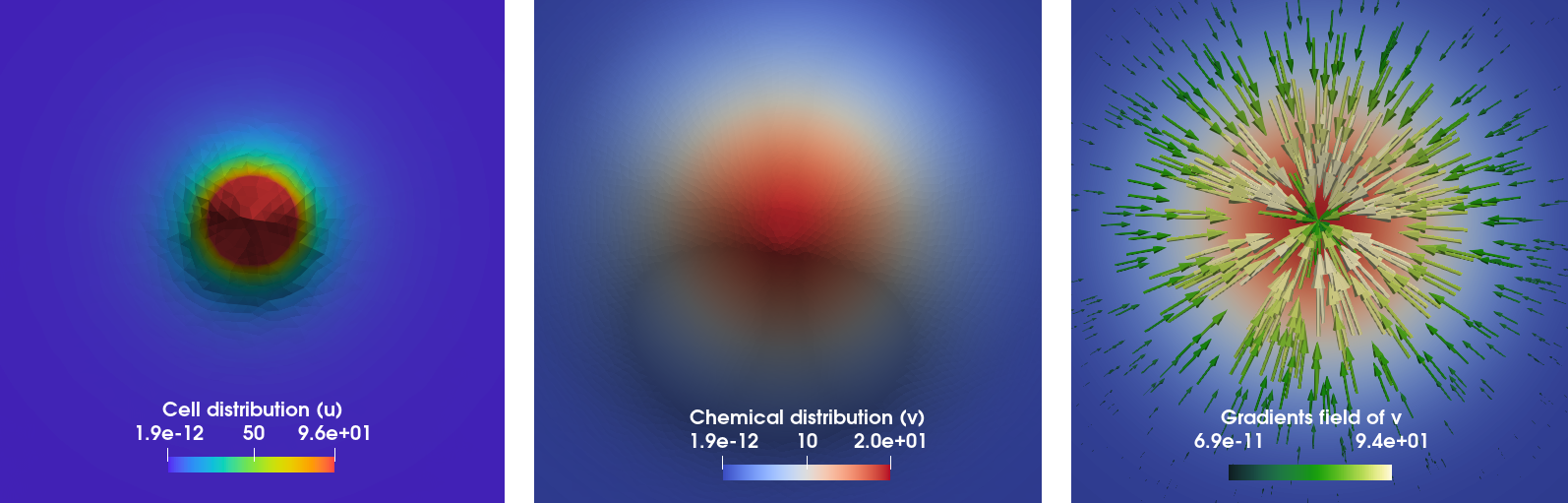}
  \\
\rotatebox{90}{\quad $t=0.0010$}&\includegraphics[width=0.73\linewidth]{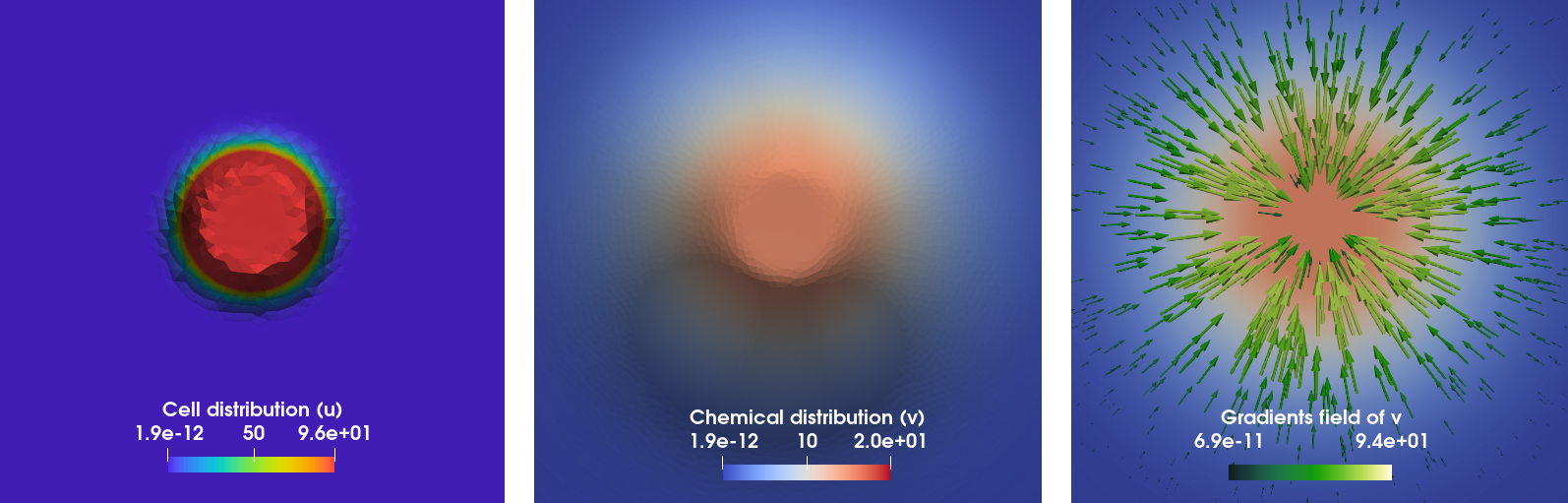}
  \\
\rotatebox{90}{\quad $t=0.0050$}&\includegraphics[width=0.73\linewidth]{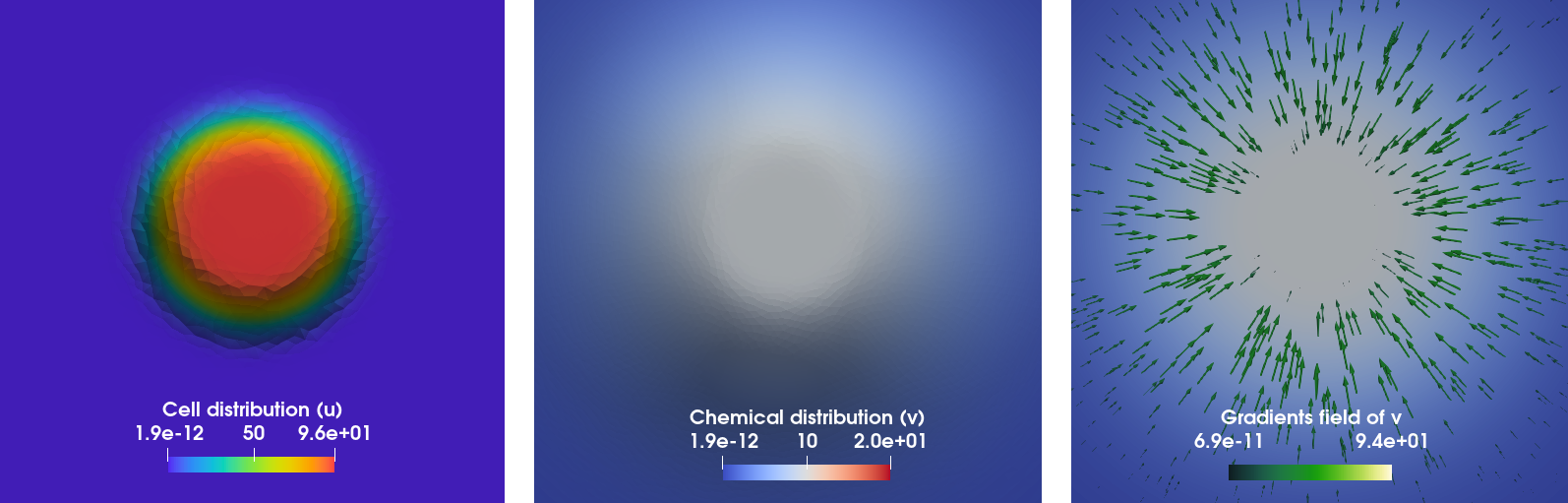}
\end{tabular}
\caption{Simulations of model \eqref{problemOriginalKSCosnumption}. Evolution of $u$ (left column), $v$ (center column) and
  $\nabla v$ (right column) in a two-dimensional domain. Each row represents their corresponding distributions, at the indicated time step. The quantitative values are deductible from each color bar.}
\label{fig:test-qualitative-2d}
\end{figure}
The first test focuses on the two-dimensional case for the classical consumption model \eqref{problemOriginalKSCosnumption}. We consider a mesh of the unit circle
$\Omega=\{(x_1,x_2) \in\mathbb{R}^2 \,\mid\, x_1^2+x_2^2 <1 \}$, meshed by
84673 triangular elements, with size $h\simeq 0.015$, and a time discretization $\Delta t=0.00001$.  We take $\chi=20$ and the Gaussian bell-shaped initial values
\begin{equation*}
u_0(x_1,x_2)=v_0(x_1,x_2)=20\,e^{-30(x_1^2+x_2^2)}.
\end{equation*}
As said, the theory predicts that no blow-up occurs in the $2d$ case; thereafter, this result may be employed to verify the consistency of our program. Indeed, in Figure~\ref{fig:test-qualitative-2d} one can see that the attraction of $v$ leads to a growth of $u$, which reaches its maximum value at $t= 0.0001$. Then, due to the consumption effect from the second equation of ~\eqref{problemOriginalKSCosnumption}, big values of $u$ imply a reduction of those of $v$ and of $\nabla v$. This ends up allowing $u$ to decrease, due to the diffusive term. This can be clearly seen again in Figure~\ref{fig:test-qualitative-2d}, where the maximum of $u$ decreases and a plateau is produced at $t=0.0010$. The decrease of $v$ and of its gradient becomes more evident for higher values of time.
\subsubsection{Three-dimensional numerical tests.}
\label{sec:3d-numerical-test}
\begin{figure}
\centering
\begin{tabular}[c]{cc}
  \includegraphics[width=0.35\linewidth]{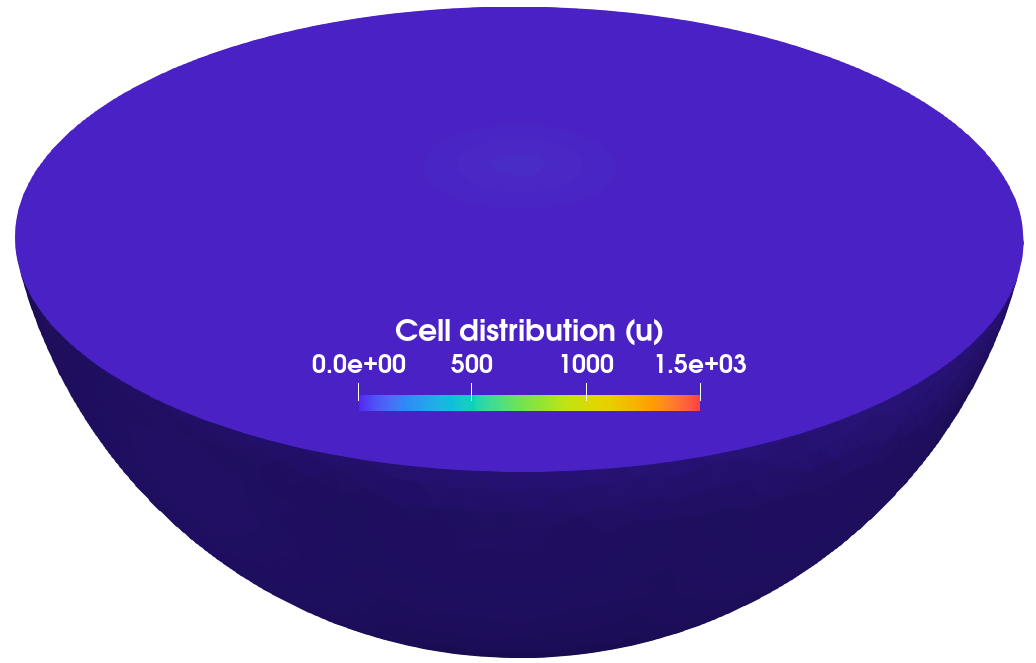}
  &\includegraphics[width=0.35\linewidth]{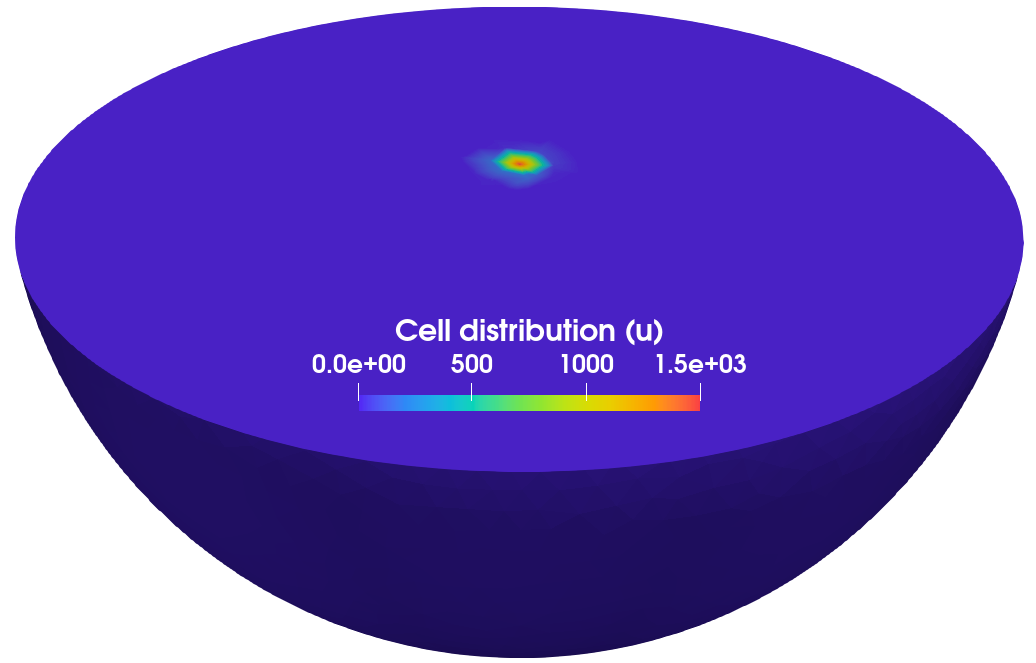}
  \\
  \includegraphics[width=0.35\linewidth]{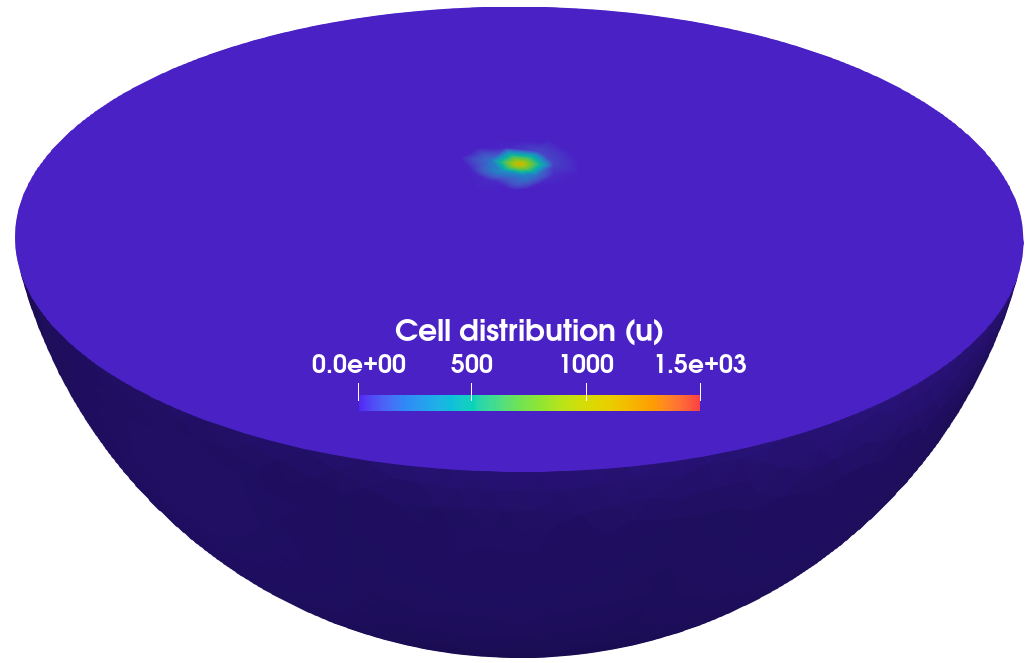}
  &\includegraphics[width=0.35\linewidth]{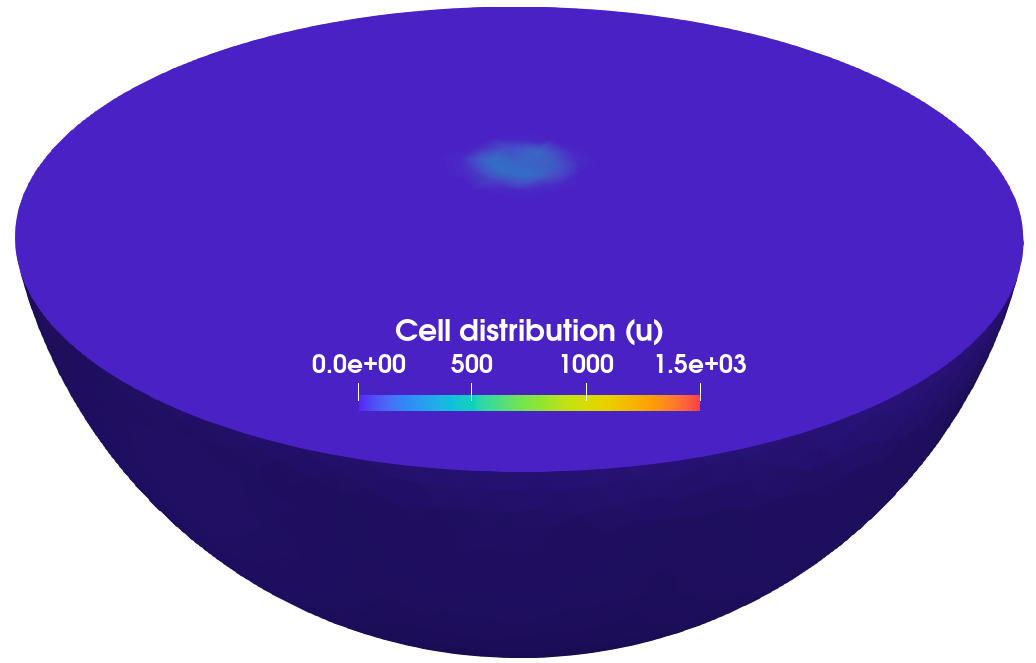}
\end{tabular}
  \caption{Simulations of model \eqref{problemOriginalKSCosnumption}. Evolution of the cells density in a three-dimensional domain: $u$ at $t=0$ (top left), $u$ at $t=0.005$ (top right), $u$ at $t=0.01$ (bottom left) and $u$ at $t=0.05$ (bottom right). The quantitative values are deductible from the color bar.}\label{fig:test-3d}
\end{figure}

Once our program is validated in a two-dimensional domain, we apply it in $3d$, where some theoretical aspects remain without a response and have to be still investigated. It is worthwhile to emphasize that even though three-dimensional simulations were highly demanding from a computational point of view, the results shown below are quite satisfactory and support our methodology. We repeat the same test described in $\S$\ref{sec:2d-numerical-test}, but on a unit sphere
$\Omega=\{(x_1,x_2,x_3) \in\mathbb{R}^3\,\mid\, x_1^2+x_2^2+x_3^2 <1 \}$, whose mesh is composed by $159296$ elements (tetrahedrons) with maximum size $h\simeq 0.11$. We take as initial data 
\begin{equation}\label{u0andv0inSimulationsThreeD}
u_0(x_1,x_2,x_3)=v_0(x_1,x_2,x_3)=20\,e^{-30\*(x_1^2+x_2^2+x_3^2)}.
\end{equation}
 Moreover, being 
$\lVert v_0\rVert_{L^\infty(\Omega)}=20$, we also fix $\chi=20 \gg 0.055536 \simeq\frac{\pi}{20 \sqrt{8} }$ (recall the condition for global boundedness $\chi \lVert v_0\rVert_{L^\infty(\Omega)}<\frac{\pi}{\sqrt{2(n+1)}}$ discussed in $\S$\ref{IntroDiscussionSection}). In this case, we have no theoretical result confirming or denying existence of unbounded solutions, but the simulations do not succumb to blow-up phenomena for model \eqref{problemOriginalKSCosnumption}  nor in $3d$. Specifically, the results are qualitatively similar to the two-dimensional ones. As a matter of fact, as contextualized in Figure~\ref{fig:test-3d}, the initial cells' density grows until reaching its maximum value at $t=0.005$; then, as the consumption effect weakens the attractive action of the chemical $v$, further growth is prevented and $u$ decreases (invalidating a possible blow-up mechanism), as specially shown at $t=0.01$ and $t=0.05$.
\subsubsection{Two- and three-dimensional numerical tests.}
\begin{figure}
  \begin{center}
    \includegraphics[width=0.8\linewidth]{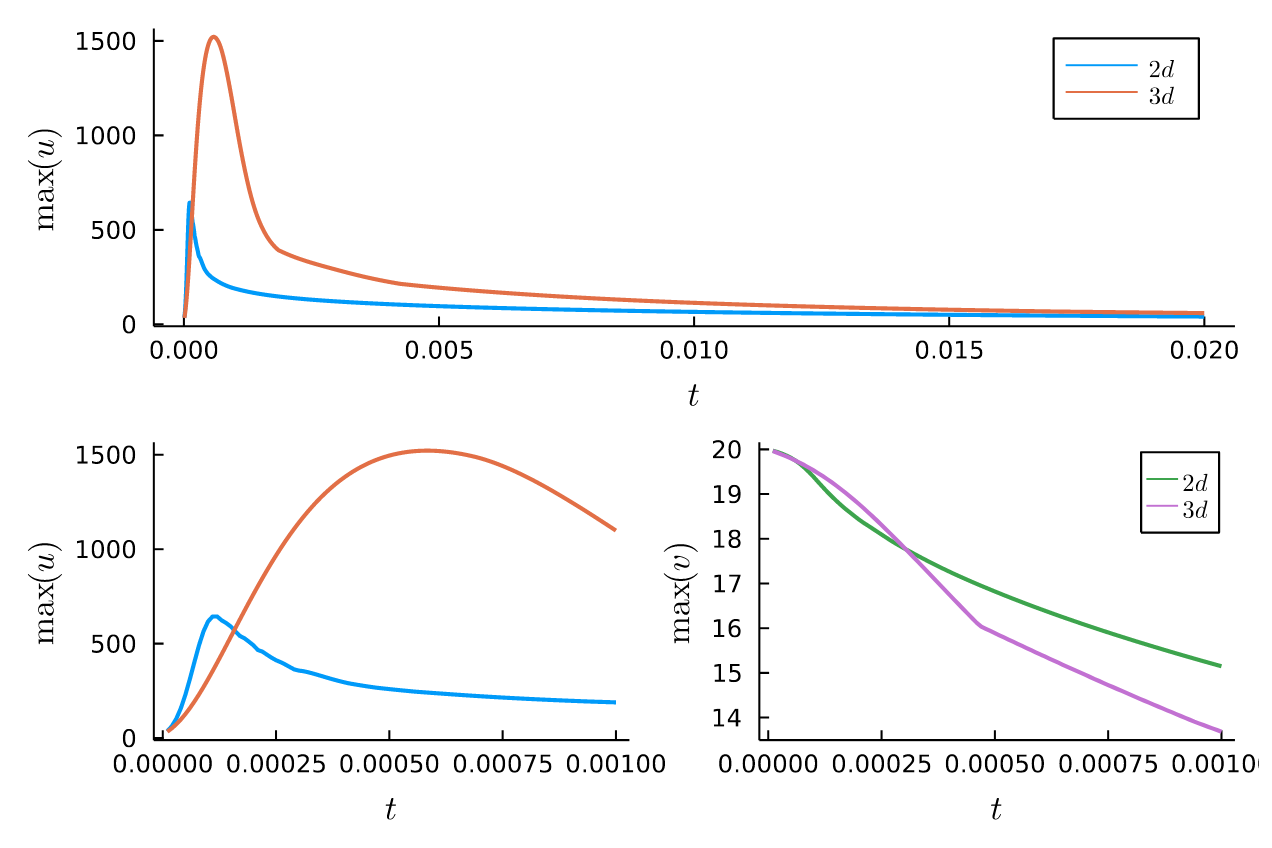}
  \end{center}
  \caption{Simulations of model \eqref{problemOriginalKSCosnumption}. Comparison of the two-dimensional and three-dimensional case along time. Evolution of the maximum of $u$ (top).  Maximum of $u$ (bottom left) and maximum of  $v$ (bottom right) at early time steps.}\label{fig:comparison-2d-3d}
\end{figure}

Figure~\ref{fig:comparison-2d-3d} indicates that in two-dimensional settings, the cells' distribution initially suffers a stronger spike than that detected in $3d$; in particular the maximum of $u$ in $2d$ is reached earlier than in $3d$. This appears consistent with the behavior of chemical agent $v$, whose maximum at the beginning is higher in $3d$. On the other hand, for the same dynamics, the maximum value of $u$ is achieved in the three-dimensional scenario; this is consistent with the well known phenomenon according to which in high dimensions blow-up phenomena in chemotaxis models are more conceivable than in lower dimensions.  Additionally, when time increases, constant steady states for $u$ and $v$ (data not shown for the chemical) are achieved. 
\subsection{Simulations of model \eqref{problem}}\label{SubSectionModelDoubleConsSimNumer}
Through the following simulations, we essentially aim at showing that the introduction in the classical consumption mechanism of a further absorbed repulsive-chemical signal $w$, enforces smoothness of solutions to the model even for large values of the chemoattractant $\chi$, associated to $v$, and a small chemorepellent $\xi$, connected to $w$.
\subsubsection{Two- and three-dimensional numerical tests.}
\begin{figure}
  \begin{center}
    \includegraphics[width=0.8\linewidth]{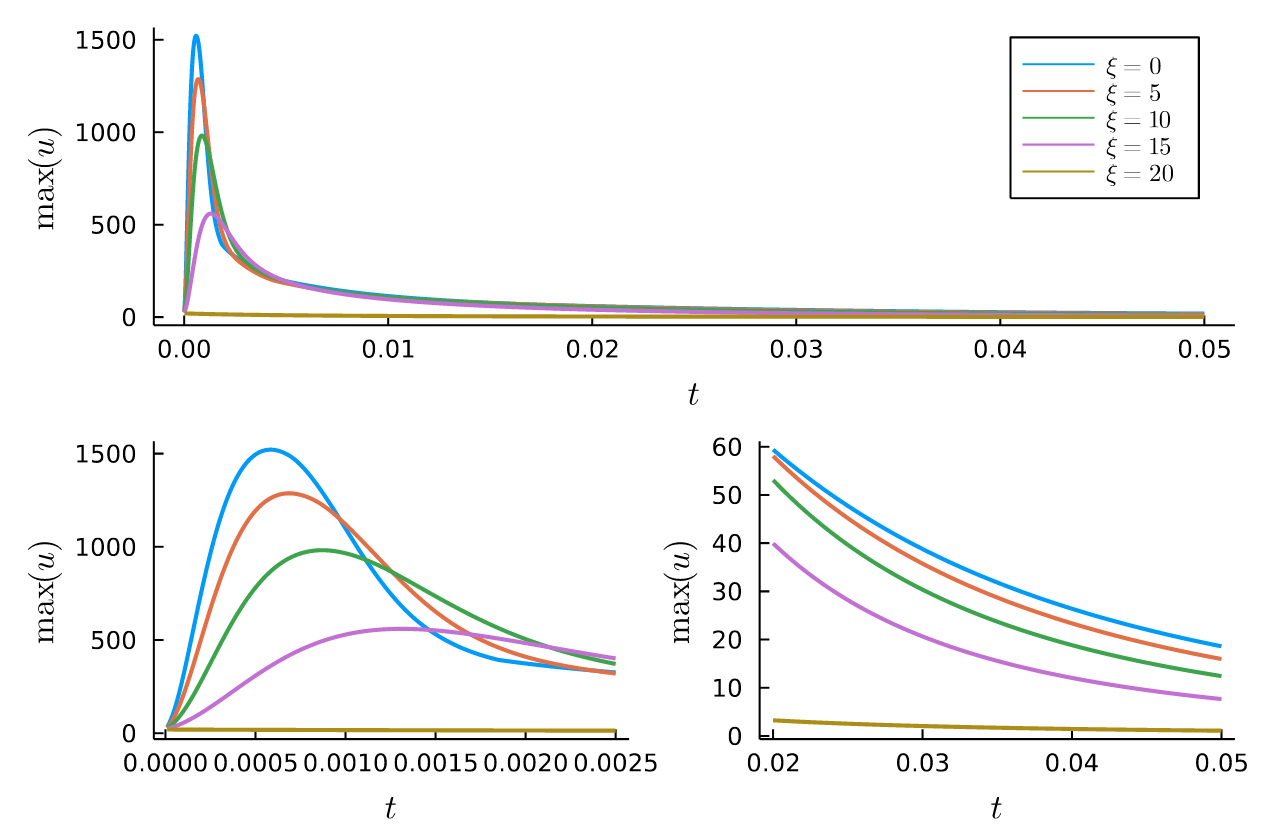}
  \end{center}
  \caption{Simulations of model \eqref{problem}. Analysis of the three-dimensional case for different values of $\xi$.}
 \label{fig:comparison-xi}
\end{figure}

In this section we discuss behaviors of the solution to system \eqref{problem} in the two-  and three-dimensional cases. Once the value of $\xi>0$ is fixed, the remaining data used in these tests are the same than those employed in $\S$\ref{sec:2d-numerical-test}, and precisely $\chi=20$ and $u_0(x_1,x_2,x_3)=v_0(x_1,x_2,x_3)=w_0(x_1,x_2,x_3)$ as in expression \eqref{u0andv0inSimulationsThreeD}. In this sense, Figure~\ref{fig:comparison-xi} allows us to discuss the evolution of the maximum of $u$ for different values of the parameter $\xi$. (Note that as much $\chi=20$ as well as any value of $\xi$ sensitively surpass the quantity $0.0033\simeq\frac{1}{15 \cdot 20}$, computed according to the assumptions in Theorem \ref{MainTheorem}.)  As expected, a simple parameter sweep shows that the maximum of $u$ is higher/lower when $\xi$ is small/large. Nevertheless, this one-to-one correspondence occurs only at the beginning and at the end of the simulations; for intermediate time steps, the curves drawing the maximum of $u$ intersect each other without any clear pattern. In particular, for (the limit case) $\xi=0$ no repulsion effect on the cells' motion takes part on the dynamics, and the most pronounced spike is produced. Conversely, when $\xi=\chi=20$, attraction and repulsion symmetrically counteract each other, so that the only  appreciated effect is reduced to the pure diffusion of living cells, with no chemotaxis phenomenon.

\subsubsection*{Acknowledgments}
SF and GV are members of the Gruppo Nazionale per l'Analisi Matematica, la Probabilit\`a e le loro Applicazioni (GNAMPA) of the Istituto Nazionale di Alta Matematica (INdAM) and are partially supported by the research project \emph{Evolutive and Stationary Partial Differential Equations with a Focus on Biomathematics}, funded by Fondazione di Sardegna (2019). GV is partially supported by MIUR (Italian Ministry of Education, University and Research) Prin 2017 \emph{Nonlinear Differential Problems via Variational, Topological and Set-valued Methods} (Grant Number: 2017AYM8XW). RRG has been supported by Grant PGC2018-098308-B-I00, by MCI N/AEI/10.13039/501100011033 and by \emph{ERDF a way of making Europe.}

\end{document}